# DISCUSSION

By Robert Serfling[1] and Yijun Zuo[2]

*University of Texas at Dallas and Michigan State University*

With delight we most heartily congratulate Hallin, Paindaveine and Šiman (HPS) on a superb and stimulating paper. It uniquely impacts our thinking about regression quantiles, multivariate quantiles, and the halfspace depth. Here we examine this highly significant contribution from the standpoints of some perspectives on multivariate quantile and depth functions, some criteria to consider in choosing such functions, and some further points about the much-studied halfspace depth. We also raise a few technical issues and questions for consideration.

**General perspectives on quantile and depth functions.** In thinking about any new contribution to multivariate quantile functions, we may draw upon the following perspectives, which also clarify the univariate case in some respects.

(P1) *In multivariate analysis, orientation to a "center" compensates for lack of a natural order.*
(P2) *In the context of quantiles, the role of "center" is naturally given to the "median."*
(P3) *The inverse of a quantile function is not the distribution $F$ but rather the rank function.*
(P4) *Depth, outlyingness, quantile, and rank functions are equivalent (DOQR paradigm).*
(P5) *Quantile functions are best viewed as parameters or characteristics of the distribution $F$.*
(P6) *Equivalence between distribution and quantile functions is not an essential requirement.*

Received August 2009.
[1]Supported by NSF Grants DMS-01-03698 and DMS-08-05786, and NSA Grant H98230-08-1-0106.
[2]Supported by NSF Grants DMS-02-34078 and DMS-05-01174.

Let us briefly elaborate on some of these points.

(P3). In the univariate case, a natural linear order makes it convenient and straightforward to define distribution and quantile functions as mutual inverses, $F$ and $F^{-1}$. However, for extension to higher dimension, the equivalent *median-oriented* formulation is the most appropriate point of departure. That is, via $u = 2p - 1$, the usual quantile function $F^{-1}(p)$, $0 < p < 1$, may be represented as $Q(u, F) = F^{-1}(\frac{1+u}{2})$, $-1 < u < 1$. Each point $x \in \mathbb{R}$ has a quantile representation $x = Q(u, F)$ for some choice of $u$, the *median* corresponding to $Q(0, F)$. For $u \neq 0$, the index $u$ indicates through its sign *the direction of $x$ from the median* and through its magnitude the *outlyingness of $x$ from the median*. Moreover, as the (unique when $F$ is continuous) solution of the equation $x = Q(u, F)$, the index $u$ defines the usual (centered) *rank function*, $R(x, F) = 2F(x) - 1$, which is thus the inverse of $Q$. It is a convenient coincidence that $R$ and $F$ are equivalent in the univariate case.

Passing to a distribution $F$ on $\mathbb{R}^d$, we may introduce associated *quantile functions* by various means, and by (P1) and (P2) it is their median-oriented formulations that are most apropos. A quantile function, indexed by $\mathbf{u}$ in the unit ball $\mathbb{B}^{d-1}(\mathbf{0})$ in $\mathbb{R}^d$, attaches to each point $\mathbf{x}$ a *quantile representation* $\mathbf{Q}(\mathbf{u}, F)$ and generates *nested* contours $\{\mathbf{Q}(\mathbf{u}, F) : \|\mathbf{u}\| = c\}$, $0 \leq c < 1$. For $\mathbf{u} = \mathbf{0}$, the most "central" point $\mathbf{Q}(\mathbf{0}, F)$ is interpreted as a $d$-*dimensional median* $\mathbf{M}_F$. Otherwise the index $\mathbf{u}$ represents a *direction* in some sense, for example, direction to $\mathbf{Q}(\mathbf{u}, F)$ from $\mathbf{M}_F$, or *expected* direction to $\mathbf{Q}(\mathbf{u}, F)$ from random $\mathbf{X} \sim F$. The magnitude $\|\mathbf{u}\|$ measures *outlyingness*, with higher values corresponding to more extreme points. The index $\mathbf{u}$ as solution of the equation $\mathbf{x} = \mathbf{Q}(\mathbf{u}, F)$ thus defines on $\mathbb{R}^d$ a *rank function* $\mathbf{R}(\mathbf{x}, F)$ which is the inverse of $\mathbf{Q}(\mathbf{u}, F)$.

(P4). From (P3) we see that quantile and rank functions, $\mathbf{Q}$ and $\mathbf{R}$, are mutually inverse, and that an *outlyingness function* is generated via $O(\mathbf{x}, F) = \|\mathbf{R}(\mathbf{x}, F)\|$. Also, an associated *depth function* $D(\mathbf{x}, F)$ measuring *centrality*, and thus inverse to $O(\mathbf{x}, F)$, is defined by some one-to-one correspondence such as $D = a + bO$ or $D = 1/(1 + O)$. Since $Q$ and $R$ are mutually inverse and $D$ and $O$ are mutually inverse, and these are linked by $O = \|R\|$, these four functions are essentially equivalent and generate the same system of contours in $\mathbb{R}^d$. Of course, the four functions have different practical roles, each with a special appeal and purpose. However, these roles are linked, and when we examine any one of $D$, $O$, $Q$ or $R$, it is important to consider it as but one element of a particular DOQR combination that is most productively viewed in its entirety. For detailed discussion and illustration with particular DOQR combinations, see [14].

Useful constructs are the associated contours, which represent equivalence classes of points of equal outlyingness (or equal depth). Using depth function contours, extensions of the univariate boxplot yield for $F$ on $\mathbb{R}^d$ notions of a "middle half" or a "middle 90%" of the population. See [7] for some general discussion. The contours do not replace, however, the underlying pointwise functions, which have their own special applications. For example, the hypothesis $H_0 : \mathbf{M}_F = \theta_0$ may be tested by the sample rank function evaluated at $\theta_0$, that is, $\mathbf{R}(\theta_0, \mathbb{X}_n)$, which represents a notion of multivariate sign test statistic. See also [7] for a great variety of nonparametric multivariate statistical methods formulated using pointwise depth functions.

(P5). The univariate median has many interesting equivalent characterizations, a number of which have been generalized to yield distinctive, meaningful notions of multidimensional median. The same is true for all the univariate quantiles. These extensions comprise a host of parameters of a multivariate $F$, each having special appeal. It is indeed useful and productive, therefore, to allow a variety of multivariate quantile functions and corresponding DOQR combinations. We may think of a quantile (or rank, depth or outlyingness) function as a "foundational" parameter, in terms of which important descriptive measures for location, spread, skewness, kurtosis, etc., may be defined. For illustration using the *spatial quantile function*, see [1] and [12]. We also mention the variety of depth-trimmed means that have been formulated and studied, for example using halfspace depth [8, 9, 11] or projection depth [21].

(P6). For the role of a quantile function $Q$ as a (sophisticated) "parameter" of $F$ as per (P5), it is not necessary or even germane that $Q$ determine $F$. For example, even though the spatial quantile function does determine $F$ [4], this aspect plays no particular role in applications. Also, there are various partial results on the degree to which the halfspace depth contours determine $F$ [5], but these have no specific role with data. In general, a "parameter" need not carry any further information about $F$ beyond that which is useful for a particular goal or purpose. Similar remarks apply to sample versions: particular statistics of interest need not retain all of the "information" in the data, and if they happen to do so, this does not guarantee that the information is organized in the most useful way.

**Some criteria for multivariate quantile and related functions.** Here we mention without elaboration some criteria which speak for themselves and arise quite typically in practical considerations. They are listed in no particular order, because their relative priorities depend upon the particular context and user.

(C1) *Equivariance of quantile and rank functions, invariance of depth and outlyingness functions.* (See [14] for elaboration.)
(C2) *Relationship between "median" and "center" relative to various notions of "symmetry."*
(C3) *Robustness.*
(C4) *Computational ease with respect to both $d$ and $n$.*

(C5) *Intuitive appeal.*
(C6) *Basis for meaningful descriptive measures for location, spread, asymmetry, kurtosis, etc.*
(C7) *Availability of applicable distribution theory, both fixed sample size and asymptotic.*
(C8) *Smoothness of contours.*
(C9) *Broadness of applicability in nonparametric sense.*

**Two new points about halfspace depth.** The quantile function constructed by HPS corresponds to the *halfspace depth*, and by the DOQR paradigm P4 the properties and behavior of the halfspace depth therefore carry over to the entire DOQR combination. The halfspace or Tukey depth has received considerable study and some of its properties and roles are alluded to by HPS as well as in the above discussion. Here we mention two further aspects that may give pause to unqualified adoption of the halfspace DOQR combination as the predominant method of choice.

(H1) *"Multivariate Tukey" $\neq$ "univariate Tukey."*
(H2) *The Tukey outlyingness function is not competitive with respect to masking breakdown point.*

Brief elaboration follows.

(H1). What is now called the halfspace depth was introduced by Tukey [15] as a method of constructing multivariate analogues of the univariate order statistics and a multivariate notion of "centrality." The corresponding outlyingness function reduces in the univariate case to a function based on tail probabilities. On the other hand, Mosteller and Tukey [10] emphasize measuring univariate outlyingness of a point $x$ by a scaled deviation, for example, $(x - \text{Median})/\text{MAD}$. This is quite different from looking at tail probabilities and its multivariate generalization turns out to be the so-called *projection outlyingness* introduced in [6] and studied in detail in [17]. The main relevance of H1 in the present context is that even Tukey did not put all his eggs in the halfspace depth basket that he invented.

(H2). In a recent study [2] of several nonparametric depth-based multivariate outlier identifiers with respect to a masking breakdown point robustness criterion, the halfspace depth was found to be singularly deficient. For classifying points as outliers or not using a chosen threshold high enough to have a low false positive rate, based on the distribution of halfspace outlyingness in a contaminated normal model, just a few outliers suffice for "masking breakdown": some arbitrarily large outliers become masked (undetected).

Here we suggest a possible explanation. For $F$ $d$-variate normal, the halfspace depth is given [3] by $D_{\text{H}}(\mathbf{x}, F) = \Phi(-\|\mathbf{x}\|)$, with $\Phi$ the univariate standard normal c.d.f. A corresponding outlyingness function designed to

take values in $[0,1]$ is $O_{\rm H}(\mathbf{x},F) = 1 - 2\Phi(-\|\mathbf{x}\|)$. It then follows [2] that for $\mathbf{X} \sim F$ the c.d.f. of $O_{\rm H}(\mathbf{X},F)$ is given by

$$F_{O_{\rm H}(\mathbf{X},F)}(\lambda) = P\left(\chi_d^2 \leq \left[\Phi^{-1}\left(\frac{1+\lambda}{2}\right)\right]^2\right).$$

Now let us consider the associated density function,

$$f_{O_{\rm H}(\mathbf{X},F)}(\lambda) = \frac{\sqrt{2\pi}(1/2)^{d/2}}{\Gamma(d/2)}\left[\Phi^{-1}\left(\frac{1+\lambda}{2}\right)\right]^{d-1}.$$

For $d = 1$ this is simply the uniform density on $[0,1]$. Unfortunately, however, for $d \geq 2$, this density has a very undesirable property: it is *monotone increasing to infinity*. Therefore, for any (high) choice of outlier threshold $\lambda$, say for false positive rate 1% in a contaminated normal model, not only will the false positives reach far into the tail where the true outliers are found, but also a typical sample will have massively many points just below the threshold. Consequently, true outliers, false positives and nonoutliers with be neighboring in rather large quantity, making it easy for outliers to become masked.

The above finding is consistent with the study of robustness of halfspace depth in [18], where it is shown that the halfspace depth of a point does not contain all the information about its relative "distance" with respect to the center of the data cloud and cannot be employed directly to identify outliers among the sample points. Indeed, outliers and points on the boundary of the convex hull may all have the same depth $1/n$. This is reflected in the low breakdown point of the pointwise halfspace depth. See [18], Example 1.

The weakness of halfspace depth with respect to robustness criteria is a serious limitation. In some applications, more robust competitors are needed.

**General views on HPS.** In the DOQR paradigm of (P3) above, every point $\mathbf{x}$ in $\mathbb{R}^d$ has a quantile representation endowing each $\mathbf{x}$ with a vector $\mathbf{u}$ in $\mathbb{B}^{d-1}(\mathbf{0})$ having a meaningful directional interpretation. Thus quantile functions $\mathbf{Q}(\mathbf{u},F)$ range pointwise through $\mathbb{R}^d$, facilitating notions of multidimensional median, directional ranks relative to $F$, depth-trimmed means, depth-trimmed scatter functionals, and a host of descriptive measures, all quite similar to well-established univariate quantile-based analysis. Conversely, starting with the direction $d$-vector $\mathbf{u}$ in $\mathbb{B}^{d-1}(\mathbf{0})$, a quantile function $\mathbf{Q}(\mathbf{u},F)$ maps directions in $\mathbb{B}^{d-1}(\mathbf{0})$ onto points $\mathbf{x}$ in $\mathbb{R}^d$.

On the other hand, in the HPS scheme, the term "quantile" is given to a $(d-1)$-dimensional regression hyperplane associated with a direction vector $\mathbf{u}$ in $\mathbb{B}^{d-1}(\mathbf{0})$, rather than to a $d$-vector in $\mathbb{R}^d$. Now we ask, is this a replacement of $\mathbf{Q}$ in the above DOQR paradigm, or is it a very interesting adjunct? The answer is found by thinking about the *contours*, that is, the upper envelopes of the HPS $\boldsymbol{\tau}$-quantile regions for fixed $\|\boldsymbol{\tau}\|$. As shown by HPS, these are simply the contours of the halfspace depth. Now, as shown in [13] and [14], a system of nested contours generates a quantile function. That is, for $D(\mathbf{x}, F)$ possessing nested contours enclosing the "median" $\mathbf{M}_F$ and bounding "central regions" of form $\{\mathbf{x} : D(\mathbf{x}, F) \geq \alpha\}$, $\alpha > 0$, the depth contours induce $\mathbf{Q}(\mathbf{u}, F)$, $\mathbf{u} \in \mathbb{B}^{d-1}(\mathbf{0})$, with each $\mathbf{x} \in \mathbb{R}^d$ given a quantile representation, as follows. For $\mathbf{x} = \mathbf{M}_F$, denote it by $\mathbf{Q}(\mathbf{0}, F)$. For $\mathbf{x} \neq \mathbf{M}_F$, denote it by $\mathbf{Q}(\mathbf{u}, F)$ with $\mathbf{u} = p\mathbf{v}$, where $p$ is the probability weight of the central region with $\mathbf{x}$ on its boundary and $\mathbf{v}$ is the unit vector toward $\mathbf{x}$ from $\mathbf{M}_F$. In this case, $\mathbf{u} = \mathbf{R}(\mathbf{x}, F)$ indicates direction toward $\mathbf{x} = \mathbf{Q}(\mathbf{u}, F)$ from $\mathbf{M}_F$, and the outlyingness parameter $\|\mathbf{u}\| = \|\mathbf{R}(\mathbf{x}, F)\|$ is the probability weight of the central region with $\mathbf{Q}(\mathbf{u}, F)$ on its boundary. All of the various depth functions considered in [7] and [19], for example, induce associated outlyingness, quantile, and rank functions. Of course, the mapping linking the two indexings $\boldsymbol{\tau}$ and $\mathbf{u}$ is not immediately transparent.

Thus the HPS quantiles induce the halfspace depth contours, which in turn induce the full halfspace DOQR combination. So, in the HPS scheme, one need not give up the usual notion of multivariate quantiles. Rather, one still may arrive at the DOQR setup for use in its intended range of applications, while at the same time enjoying additional benefits provided by the hyperplanes. These may be regarded as an adjunct to the DOQR paradigm in the halfspace case. It then becomes of interest to explore the possibility of such adjuncts relative to other choices of depth function.

Salient features of the HPS quantile approach, in terms of (P1)–(P6) and (C1)–(C9), are thus evaluated in terms of the halfspace depth. On balance, the halfspace depth is among a handful of leading depth functions, of which no single one predominates, the priorities among different perspectives and criteria depending on the context and the user. Overall, the halfspace depth is strong relative to (C1)–(C9) except for some limitations with respect to (C3), (C4) (but see below) and (C8).

From the standpoint of depth functions, a key contribution of HPS is to strengthen the appeal of the halfspace approach by providing it with two very important gains in its assets:

- Relative to criterion (C4), *a significant new computational pathway to halfspace depth contours.*
- Relative to criterion (C5), *a significant linkage with multi-output regression quantiles.*

*Some miscellaneous issues and questions.* We augment the preceding general views and comments with some specific issues and questions.

- *Connections with univariate quantiles and quantile analysis.* In the univariate case, the HPS lower quantile hyperplane reduces to $(-\infty, F^{-1}(\tau)) \cup$ $(F^{-1}(1-\tau),+\infty)$, the complement of the "$\tau$ depth contour." This set comprises the upper and lower tails of probability $\tau$, with total probabiility $2\tau$, consistent with the conventional univariate case. However, the discussion of equation (3.2a) seems to indicate that this probability should be just $\tau$, seemingly a contradiction.

The discussion of the univariate case following Definition 2.1 is a bit sketchy, using vector $\boldsymbol{\tau}$ notation instead of scalar $\tau$ and also throwing in a rather assertive statement, without qualification or elaboration, that depth contours should be associated with contour-valued rather than point-valued quantiles. Regarding the latter, we have shown in our discussion of (P3) how these contours do not dispense with the points that comprise them, and it would be awkward to insist that all of univariate quantile usage be revised to use only contours and avoid point-wise quantiles. We endorse keeping everything in sight, both in the univariate case and in multivariate extensions.

Equation (3.9) imposes a restriction on $\tau$ that makes empirical versions well-defined, but this leaves them evidently undefined for $\tau$ values in $(0, N/n)$ and $(1-P/n, 1)$. However, classical univariate empirical quantiles are well defined for all $\tau \in (0,1)$.

- *Moment and regularity assumptions.* We note that the regularity and moment conditions imposed by Assumptions (A) and $(A'_n)$ are uncompetitively strong for a quantile and depth function methodology. After all, univariate quantile analysis requires neither regularity nor moment assumptions to get started, and such results as Bahadur representations for sample quantiles require second-order regularity but not moment assumptions. Likewise, halfspace depth is well-defined without moment assumptions. Thus Assumptions (A) and $(A'_n)$ represent an additional price to be paid for the hyperplane quantiles and behavior of sample versions.
- *The assumption of "general position."* Since depth contours are well-defined for any data set (not necessarily in general position), we query whether this assumption is necessary in Theorem 4.2. Also, in that theorem, $\ell/n$ must be less than the maximum halfspace depth (see [3] for related discussion on maximum possible halfspace depth).
- *Compactness of the $R(\tau)$ regions.* The discussion following Theorem 4.2 includes the statement that the supremum of all $\tau$ such that $R(\tau) \neq \varnothing$ belongs to the interval $[1/(k+1), 1/2]$. Here $1/(k+1)$ should be replaced by $1/n$, since in $\mathbb{R}^k$ we may suppose that the data points are in general position and on the vertices (corners) of the hyperpolygon, and then for halfspace depth the supremum of all $\tau$ is $1/n < 1/(k+1)$.

We note in passing that for halfspace depth the characterization that this supremum is $1/2$ if and only if the distribution of $\mathbf{Z}$ is angularly symmetric has been treated earlier in detail in [16] and [20].

- *Potential practical applications of the hyperplane quantiles.* Applications of halfspace depth, its contours, and its related Q, O, R functions, are quite well established and familiar. It is of interest to know how the hyperplane quantiles and their sample versions would be used in practice. What added methodology is acquired using these particular entities? And if existing pointwise depth and quantile methods are to be de-emphasized or reformulated, how do the hyperplane methods accomplish all the same goals?

**Summary.** This paper extends regression quantiles to the multivariate setting and links with the well-established halfspace depth. With respect to the latter, significant new insights and computational approaches are provided. The paper treats its topic with great thoroughness and flair and, indeed, is a *tour de force*.

In a larger view, the DOQR paradigm in the halfspace case is augmented by an additional entity, a directional hyperplane. It is of interest to know more about practical applications of the "H" in this extended "DOQRH" paradigm, and it may also be of interest to explore whether a "DOQRH" paradigm has meaningful formulation in the context of other depth functions.

The research community truly owes Marc, Davy and Miroslav a great debt of gratitude for their outstanding work that changes our perceptions, adds to our tools, and stimulates interesting further inquiries. We look forward to further developments!

Department of Mathematical Sciences
University of Texas at Dallas
Richardson, Texas 75083-0688
USA
E-mail: serfling@utdallas.edu

Department of Statistics and Probability
Michigan State University
East Lansing, Michigan 48824
USA
E-mail: zuo@msu.edu